\documentclass{amsart}
\usepackage{latexsym,amsmath,amssymb}
\newtheorem{thm}{Theorem}[section]

\newtheorem{lem}[thm]{Lemma}
\newtheorem{exam}[thm]{Example}
\newtheorem{prop}[thm]{Proposition}
\theoremstyle{definition}\newtheorem{definition}[thm]{Definition}
\theoremstyle{remark}
\newtheorem{rem}[thm]{Remark}
\numberwithin{equation}{section}

\begin{document}

\title[]
{FIXED POINT OF SUBADDITIVE MAPS and some non-linear integral equations  }

\author{\sc Yousef Estaremi and Bahman Moeini }
\address{\sc yousef estaremi }
\email{yestaremi@pnu.ac.ir} \email{}
\address{Department of mathematical
sciences, Payame Noor University, Tehran, Iran}
\address{\sc bahman moeini}
\email{moeini$_{-}$summer@kiau.ac.ir}
\address{Department of mathematics, Hidaj Branch, Islamic Azad University,
Hidaj, Iran}

\thanks{}

\thanks{}



\date{}


\commby{}


\begin{abstract}
In this paper,  first some results of \cite{be3} are extended for
subadditive separating maps between $C(X,E)$ and $C(Y,E)$, such
that $E$ is a unital Banach algebra. Then we give some conditions under which a
strongly subadditive map has a unique fixed point. Finally as an application
the existence and uniqueness of solution for a nonlinear integral equation is discussed.

 \noindent {}
\end{abstract}

\maketitle
\section {\sc Introduction}
Let $X$ and $Y$ be compact Hausdorff spaces and let $E$ be a
unital Banach algebra. $C(X,E)$ and $C(Y,E)$ denote the spaces of
 $E$-valued continuous functions on $X$ and $Y$,
respectively. A map $H:C(X,E)\rightarrow C(Y,E)$ is separating if
$\|f(x)\|.\|g(x)\|=0$  implies that $\|Hf(y)\|.\|Hg(y)\|=0$ for
all  $f,g\in C(X,E)$  and  $x\in X, y\in Y$. The study of separating maps between different spaces
of functions(as well as operator algebras) has attracted a
considerable interest in resent years, see for example
\cite{aj,be2,ja,vv} and references therein.
The well known results concerning separating maps from $C(X)$ to $C(Y)$, for compact Hausdorff
spaces $X$ and $Y$, have been extended to not necessarily linear case in \cite{be3,se,vv}.
In \cite{se} we considered the results of \cite{be3} and extended some of
them to the case that $H$ is a subadditive separating map between regular Banach function algebras
$A$ and $B$ on $X$ and $Y$, respectively.
Also, in this paper we are going to investigate some properties of subadditive
separating maps between $C(X,E)$ and $C(Y,E)$, where
$E$ is a unital Banach algebra. Mainly we will give some conditions under which the strongly
subadditive map $H$ has a unique fixed point. As an application, we use the fixed point to prove
the existence of unique solution of a Volterra type integral equation.
\section {\sc  SUBADDITIVE SEPARATING MAPS BETWEEN $C(X,E)$ and $C(Y,E)$}
Throughout this paper $X$ and $Y$ are compact Hausdorff spaces and
$E$ is a unital Banach algebra with unit element $e_0$. Given an
element $f\in C(X,E)$, let $coz(f)$ denote the cozero set of $f$
i.e.,  $coz(f)=\{x\in X:f(x)\neq0\}$,  $ClU$ denotes the topological closure of the set $U$ and $CU$
denotes the complement of the set $U$. For a separating map $H:C(X,E)\rightarrow C(Y,E)$ and $y\in Y$,
$\hat{y}oH:C(X,E)\rightarrow E$ and $\hat{y}oH(f)=Hf(y)$, for all
$f\in C(X,E)$. Here we recall some concepts that were defined in \cite{be3}.
\begin{definition}\label{d2.1}\cite{be3},  a) A map
$H:C(X,E)\rightarrow C(Y,E)$ is called subadditive, if  $$\|H(f+g)\|\leq\|Hf\|+\|Hg\|,$$
for all $f,g\in C(X,E)$.\\
b) Let $M$ be positive and $H$ be subadditive map. If for each
$f\in C(X,E)$ there exists $\varepsilon>0$ such that
$\|Hf-Hg\|\leq M\|H(f-g)\|$,
 for all $g\in C(X,E)$ satisfying $\|f-g\|<\varepsilon$, then $H$ is called strongly subadditive.
\end{definition}
\begin{exam} Let $g:Y\rightarrow X$ be continuous and let $w\in C(Y,\mathbb{C})$. The map
\[
H:C(X,\mathbb{C})\rightarrow C(Y,\mathbb{C}), \  \  \ \ \  f\rightarrow |w.(f\circ g)|
\]
 is separating and strongly subadditve, but not additive.
\end{exam}
\begin{definition}\label{d2.3}\cite{be3}, A map
$H:C(X,E)\rightarrow C(Y,E)$ is called pointwise subadditive if
\[
\|H(f+g)(y)\|\leq\|Hf(y)\|+\|Hg(y)\|,
\]
for all $f,g\in C(X,E)$ and $y\in Y$,
\end{definition}
\begin{definition}\label{d2.4}\cite{be3}, Let
$H:C(X,E)\rightarrow C(Y,E)$ be a separating map. An open subset U
of $X$ is called a vanishing set for $\hat{y}oH$ if for each $f\in
C(X,E)$, $coz(f)\subseteq U$ implies that  $\|\hat{y}H(f)\|=0$.
The support of $\hat{y}oH$ is then defined by
\begin{center}
supp $\hat{y}oH=X\backslash\{V\subseteq X:V$ is a
vanishing set for $\hat{y}oH\}$.
\end{center}
\end{definition}
For $e\in E$ we define the function $1_{e}:X\rightarrow E$, such
that $1_{e}(x)=e$ for all $x\in X$. Clearly $1_{e}\in C(X,E)$. In
this section we assume that $\|H1_{e_0}(y)\|\neq0$ for all $y\in Y$.
Since $\|1_{e_0}\|.0=0$ and $H$ is separating, then $H0=0$.
The definition of functions $1_{e}$ is the main tools in this section
to extend scalar-valued case to vector valued case.
\begin{rem}\label{r2.5} (Decomposition of the
identity) If $K=\mathbb{R}$ or $\mathbb{C}$, note that for any finite
cover $\{U_{i}\}^{n}_{i=1}$ of open subsets of $X$,
there is a continuous decomposition of identity, $\{e_{i}\}^{n}_{i=1}\subseteq C(X,E)$,
subordinate to the $U_{i}$, namely $\Sigma^{n}_{i=1}e_{i}=1$ and
$coz(e_{i})\subseteq U_{i}$, for each $i=1,2, ..., n$\cite{nai1}.
\end{rem}
By using Remark \ref{r2.5} we have the next theorem.
\begin{thm}\label{t2.6}
 Let $H:C(X,E)\rightarrow
C(Y,E)$ be  pointwise subadditive separating map. Then for every
$y\in Y$, the set supp$(\hat{y}oH)$ is singleton.
\end{thm}
\begin{proof}
In the first we assume that supp$(\hat{y}oH)$ is
empty, then $X=\cup_{\alpha}U_{\alpha}$ where
$\{U_{\alpha}\}_{\alpha}$ is the collection of vanishing sets for
$\hat{y}oH$(we note that, since $H0=0$ the empty set is clearly a
vanishing set). Since $X$ is compact, then we can choose finite
number of vanishing sets $U_{1}, U_{2},...,U_{n}$, such that
$X=\cup^{n}_{i=1}U_{i}$, by using Remark \ref{r2.5}, there exist $f_{1},
f_{2},...,f_{n}\in C(X)$ with $coz(f_{i})\subseteq U_{i}$ and
$\Sigma^{n}_{i=1}f_{i}=1$ on $X$. Let $f_{i,e_{0}}:X\rightarrow
E$; $f_{i,e_{0}}(x)=f_{i}(x)e_{0}$. It is clear that
$f_{i,e_{0}}\in C(X,E)$ and $coz(f_{i,e_{0}})\subseteq U_{i}$.
Also $f=\Sigma^{n}_{i=1}f_{i,e_{0}}f$, for every $f\in
C(X,E)$.
Note that $f_{i,e_{0}}f\in C(X,E)$ and $coz(ff_{i,e_{0}})\subseteq
U_{i}$, so $H(f_{i,e_{0}}f)(y)=0$. By pointwise subadditivity of
$H$,
\begin{align*}
\|H1_{e_0}(y)\|&=\|H(\Sigma^{n}_{i=1}f_{i,e_{0}}1_{e_0})(y)\|\\
&\leq\Sigma^{n}_{i=1}\|H(f_{i,e_{0}}1_{e_0})(y)\|=0.
\end{align*}
But this is a contradiction. Thus $X\neq\cup_{\alpha}U_{\alpha}$
i.e. $X\backslash\cup_{\alpha}U_{\alpha}\neq\emptyset$. Suppose
that $x$ and $x'$ are disjoint elements of
$X\backslash\cup_{\alpha}U_{\alpha}$ and $U$ and $V$ are disjoint
open neighborhoods of $x$ and $x'$ in $X$, respectively. Since
$x\notin\cup_{\alpha}U_{\alpha}$, then $U$ can't be a vanishing
set for $\hat{y}oH$. Therefore there must exist some $f\in C(X,E)$
such that $coz(f)\subseteq U$ and $\|Hf(y)\|\neq0$. Similarly
there must exist some $g\in C(X,E)$ such that $coz(g)\subseteq V$
and $\|Hg(y)\|\neq0$. Consequently, $\|Hf(y)\|.\|Hg(y)\|\neq0$,
this contradicts the fact that $H$ is separating. This implies
that supp$(\hat{y}oH)$ is singleton.
\end{proof}
Under the hypothesis of the Theorem \ref{t2.6} we can correspond to each
$y\in Y$ an element $h(y)\in X$, which is the unique point of
supp$(\hat{y}oH)$. We call the map $h:Y\rightarrow X$, defined in
this way, support map of $H$.
The next theorem can also be obtained with a minor modifications of the proof of
Theorem 4.3 of \cite{be3}, so we omit it's proof.
\begin{thm}\label{t2.7} Let $H:C(X,E)\rightarrow C(Y,E)$
be  pointwise subadditive separating map. Then
a) $h(coz(Hf))\subseteq suppf=Clcozf$,for all $f\in C(X,E)$.
\vspace*{0.3cm}
b) $\{h(y)\}=\cap_{\|\hat{y}oH(f)\|\neq0}suppf =supp(\hat{y}oH)$.
\end{thm}
Before investigating more properties of $h$, we define the
following concept which is introduced in \cite{be3} for a special
case. This concept will plays a role in automatic continuity
results (Theorem \ref{t2.11}).
\begin{definition}\label{d2.8} Let
$H:C(X,E)\rightarrow C(Y,E)$
be  pointwise subadditive separating map with support map h. We
say that H is strongly pointwise subadditive, if for  each $y\in
Y$ there exists $M_{y}>0$ and for each element $c \in E$ there
exists $\delta_{c,y}>0$ (depending on c and y) such that
$$\|Hf(y)-Hg(y)\|\leq M_{y}\|H(f-g)(y)\|,$$ holds for all $f,g\in
C(X,E)$ with $f(h(y))=c$ and $\|f(h(y))-g(h(y))\|<\delta_{c,y}$.
\end{definition}
For some significant example of pointwise and strongly pointwise
subadditive separating maps we refer to \cite{be3}. Pointwise strong
subadditivity is very similar to strong subadditivity; under certain conditions,
Pointwise strong subadditivity implies strong subadditivity. We leave simple proof
of the assertion in next proposition to the reader.
\begin{prop}\label{p2.9} Let $H:C(X,E)\rightarrow
C(Y,E)$ be strongly pointwise subadditive. If for each $f\in
C(X,E)$ the set $\{\delta_{f(h(y)),y}:y\in Y\}$ has an upper bound
and the set $\{M_{y}; y\in Y\}$ is bounded above as well, then H
is strongly subadditive.
\end{prop}
Here we recall the definition of detaching maps.
\begin{definition}\label{d2.10} The map $H:C(X,E)\rightarrow C(Y,E)$
is detaching, if for any two distinct point's $y, \acute{y}\in Y$,
there exist $f,g\in C(X,E)$ such that $cozf\cap cozg=\emptyset$
and $\|Hf(y)\|\|Hg(\acute{y})\|\neq0$.
\end{definition}
We show in Theorem \ref{t2.11}(b) that if $H$ is strongly pointwise
subadditive, then $h$ is continuous. Thus a separating connection
between $C(X,E)$ and $C(Y,E)$ establishes a continuous connection
between $X$ and $Y$.
\begin{thm}\label{t2.11} Let $H:C(X,E)\rightarrow
C(Y,E)$ be strongly pointwise subadditive map. Then the followings hold.\\
a) For any $f,g\in C(X,E)$ if $f=g$ on the open set
$U\subseteq X$, then $Hf=Hg$ on $h^{-1}(U)$.\\
b) The support map $h:Y\rightarrow X$ is continuous.\\
c) If $H$ is injective, then $h$ is surjective.\\
d) If $H$ is detaching if and only if $h$ is injective.
\end{thm}
\begin{proof}  a) Suppose that $f\in C(X,E)$`  and
$y\in Y$ such that $f=0$ on the nonempty open set $U\subseteq X$
where $h(y)\in U$. For each $x\in U^{c}$ there exists a vanishing
set $U_{x}$ for $\hat{y}oH$ which contains $x$. Since $X$ is
compact, then $U^{c}\subseteq \cup^{n}_{i=1}U_{i}$, where
$U_{i}$'s are vanishing sets for $\hat{y}oH$. Thus $h(y)\notin
\cup^{n}_{i=1}U_{i}$. If we set $U_{n+1}=U$, then
$X\subseteq\cup^{n+1}_{i=1}U_{i}$. Thus by Remark \ref{r2.5} and the
method that we used in the proof of Theorem \ref{t2.6}, we will have
$e_{1}, e_{2},..., e_{n+1}\in C(X,E)$ such that $coze_{i}\subseteq
U_{i}$ and $\Sigma^{n+1}_{i=1}e_{i}=1_{e_{0}}$ on $X$ and so
$\Sigma^{n+1}_{i=1}fe_{i}=f$. Since $f=0$ on $U_{n+1}$ and
$coze_{n+1}\subseteq U_{n+1}$ it follows $fe_{n+1}=0$. Also,
$coz(fe_{i})\subseteq U_{i}$, $1\leq i\leq n$. Thus we have
$H(fe_{i})(y)=0$ for $1\leq i\leq n+1$. Since H is pointwise
subadditive, then
\begin{align*}
\|Hf(y)\|&=\|H(\Sigma^{n+1}_{i=1}fe_{i})(y)\|\\
&\leq\Sigma^{n+1}_{i=1}\|H(fe_{i})(y)\|=0.
\end{align*}
Thus $Hf=0$ on $h^{-1}(U)$. Now suppose that $f,g\in C(X,E)$ and $f=g$ on U,
then $H(f-g)=0$ on $h^{-1}(U)$. By Definition \ref{d2.8}, there exist $M_{y}>0$ and
$\delta_{f(h(y)),y}>0$ for $y\in h^{-1}(U)$ such that
\[
\|Hf(y)-Hg(y)\|\leq M_{y}\|H(f-g)(y)\|=0,
\]
since $\|f(h(y))-g(h(y))\|=0<\delta_{f(h(y)),y}$. It follows that $Hf=Hg$ on
$h^{-1}(U)$.\\
b) Suppose that  $y\in Y$ and V is an open neighborhood of $h(y)$
in $X$. Compactness of $X$ implies that, there exist open
neighborhoods U of $CV$ and W of $h(y)$ such that $ClU\cap
ClW=\emptyset$. We can choose $f\in C(X,E)$ such that
$f=1_{e_{0}}$ on $ClW$ and $f=0$ on
 $ClU$. Since $f=1_{e_{0}}$ on $ClW$, then $Hf=H1_{e_{0}}$ on $h^{-1}(W)$ by (a).
Hence $Hf(y)=H1_{e_{0}}(y)\neq0$ i.e $y\in cos Hf$ and $coz(H(f))$ is a
open neighborhood of $y$. Suppose that $y'\in Y$ such that
$y'\notin h^{-1}(V)$. Since $f=0$ on $U$, it follows that from
(a), $Hf=H0$ on $h^{-1}(U)$ and so $Hf(\acute{y'})=0$. Thus
$cozHf\cap h^{-1}(CV)=\emptyset$ and so $h(cozHf)\subseteq V$. It follows that h is
continuous at y.\\
C) Suppose that $U\subseteq X$ is an open and nonempty
set. By Uryson's lemma there exist nonzero function $f\in C(X,E)$
such that $cozf\subseteq U$, since H is injective, then $Hf\neq0$
i.e $Hf(y)\neq0$ for some $y\in Y$. Thus $h(y)\in suppf\subseteq
U$ i.e $h(Y)\cap U\neq\emptyset$, it follows that $h(Y)$ is dense
in $X$. Also, $h(Y)$ is closed, because $h$ is continuous and $Y$
is compact. This implies that $h$ is surjective.\\
d) Suppose that H is detaching and $y,\acute{y}$ are distinct
point of $Y$. Thus, there exist $f,g\in C(X,E)$ such that
$cozf\cap cozg=\emptyset$ and $\|Hf(y)\|.\|Hg(\acute{y})\|\neq0$.
By Theorem \ref{t2.7} it follows that  $h(y)\neq h(\acute{y})$  i.e h is injective.
Conversely, suppose that $h$ is injective and $y\neq y'$. So
$h(y)\neq h(y')$. Choose open neighborhoods $U$ and $V$ of $h(y)$
and $h(y')$, respectively, such that $ClU'\subseteq U$ and
$ClV'\subseteq V$, where $U'$ and $V'$ are also neighborhoods of
$h(y)$ and $h(y')$, respectively. Thus we can find $f, g \in C(X,
E)$ such that $f=1_{e_0}$ on $U'$ and $f=o$ on $CU$, also
$g=1_{e_0}$ on $V'$ and $g=o$ on $CV$. By (a) we have
$\|Hf(y)\|\neq 0$ and $\|Hg(y')\|\neq 0$. This completes proof.
\end{proof}
\begin{thm}\label{t2.12} If $H:C(X,E)\rightarrow C(Y,E)$
is biseparating bijection, such that $H$ and $H^{-1}$ are strongly
pointwise subadditive, then the support map $h$ of $H$ is a
surjective homeomorphism, and the support map $H^{-1}$ is
$h^{-1}$.
\end{thm}
\begin{proof}First we show that $h$ is 1-1. If $h$
isn't 1-1, then there exist $y, y'\in Y$ such that $y'\neq y$ and
$h(y)=h(y')$. Bijectivity of $H$ implies that, there exist $f,g
\in C(X, E)$ such that $Cl(cozHf)\cap Cl(cozHg)=\emptyset$,
$Hf(y)\neq0$ and $Hg(y')\neq0$. As $H$ is biseparating, $cozf\cap
cozg=\emptyset$. By Theorem \ref{t2.7}(a), $h(y)=h(y')\in cozf\cap cozg$.
Let $w$ be the support map of $H^{-1}$. Then by Theorem \ref{t2.7}(a) we
have $w(coz H^{-1}Hf)=w(cozf)\subseteq Cl cozHf$ and $w(coz
H^{-1}Hg)=w(cozg)\subseteq Cl cozHg$. Since $w$ is continuous and
$h(y)=h(y')\in Clcozf\cap Clcozg$, then $w(h(y))=w(h(y'))\in
ClcozHf\cap ClcozHg$.
This is a contradiction. Therefore, $h$ is 1-1.\\
Since $H$ is bijective and $X$ and $Y$ are compact,
then by Theorem \ref{t2.11}, $h$ is a homeomorphism from $Y$ onto $X$. In
the end, if $Hf(y)\neq0$ then $h(y)\in Clcozf$. Since $w(h(y))\in
ClcozHf$ for all $f$ such that $Hf(y)\neq0$. But as $H$ is
bijective, $\cap_{\|Hf(y)\|\neq0}ClcozHf=\{y\}$. This implies that
$w(h(y))=y$, i.e., $w=h^{-1}$.
\end{proof}
In the sequel we consider continuity and form of
strongly pointwise subadditive separating maps. First we give some definitions.
\begin{definition}\label{d2.13} Let $E$ and $F$ be normed 
linear spaces and  $A$ be a map from $E$ into $F$. If there exists $D>0$
such that $\|Ae\|\leq D \|e\|$ for all $e\in E$, then $A$ is
norm-bounded.
\end{definition}
\begin{definition}\label{d2.14} Let $H$ map $C(X,E)$ into
$C(Y,E)$. If there exists $D>0$ such that $\|H(1_e)\|\leq
D\|A(1_{e_0})\|.\|e\|$ for all $e\in E$. Then $H$ is
$1_{e_0}$-bounded.
\end{definition}
We refer to \cite{be3} for some example of maps that are
continuous but neither $1_{e_0}$- or norm-bounded. If $H:C(X,E)\rightarrow C(Y,E)$
is norm-bounded and strongly subadditive, then by last definitions we conclude that
$H$ is continuous. Also, if $H$ is $1_{e_0}$-bounded and strongly subadditive on
the linear span [$1_{e_0}$] of $1_{e_0}$ i.e.,
$\{\sum^{n}_{i=1}e_i.1_{e_{0}} |  e_i\in E, n\in \mathbb{N}\},$
then $H$ is continuous on [$1_{e_0}$].
\begin{definition}\label{d2.15} Let $p:Y\rightarrow X$ be
continuous. For each $y\in Y$ and $f\in C(X,E)$, define
$H:C(X,E)\rightarrow C(Y,E)$ to be $Hf(y)=H(1_{f(p(y))})(y)$. $H$
is called a composition map.
\end{definition}
\begin{thm}\label{t2.16} Let $H:C(X,E)\rightarrow
C(Y,E)$ be a composition map. If $H$ is $1_{e_0}$-bounded, then
$H$ is norm-bounded.
\end{thm}
\begin{proof} We leave the proof to the reader.
\end{proof}
\begin{definition}\label{d2.17} Let $H:C(X,E)\rightarrow
C(Y,E)$ be a map. The set $Y_c=\{y\in Y: \hat{y}\circ H$ is
continuous$\}$ is called the continuity set of $H$.
Finally we explain the form of strongly point-wise subadditive separating maps on
vector valued continuous functions.
\end{definition}
\begin{thm}\label{t2.18} Let $H:C(X,E)\rightarrow
C(Y,E)$ be separating and strongly point wise subadditive with
support map $h$. Then the followings hold.\\
(a) If $y\in Y_c$, then $Hf(y)=H(1_{f(h(y))})(y)$ for all $f\in
C(X,E)$.
(b) If $H$ is $1_{e_0}$-bounded, then $y\in Y_c$ if and only if
\[
Hf(y)=H(1_{f(h(y))})(y)
\]
for all $f\in C(X,E)$.
\end{thm}
\begin{proof} Let $y\in Y_c$, $f\in C(X,E)$
and $\varepsilon>0$. Since $f$ is continuous, then there exists
the neighborhood $U_{\varepsilon}$ of $h(y)$ such that
$\|f(z)-f(h(y))\|<\varepsilon$ for all $z\in U_{\varepsilon}$. Let
$W$ be a neighborhood of $h(y)$ such that $ClW\subset
U_{\varepsilon}$. Suppose that $k_{\varepsilon}\in C(X,E)$ such
that $k_{\varepsilon}=1_{e_0}$ on $W$, $k_{\varepsilon}=0$ on
$CU_{\varepsilon}$ and $0\leq k_{\varepsilon} \leq 1_{e_0}$. As
$\varepsilon\rightarrow 0$,
$k_{\varepsilon}(f-1_{f(h(y))})\rightarrow 0$; consequently
$H(k_{\varepsilon}(f-1_{f(h(y))}))(y)\rightarrow 0$. Since $H$ is
strongly pointwise subadditivity and
$(k_{\varepsilon}f)(h(y))=(k_{\varepsilon}1_{f(h(y))})(h(y))=f(h(y))$
for all $\varepsilon>0$, then there exists $M_y>0$ such that
\[\|H(k_{\varepsilon}f)(y)-H(k_{\varepsilon}1_{f(h(y))})(y)\|\leq
M_y\|H(k_{\varepsilon}(f-1_{f(h(y))}))(y)\|\rightarrow 0,
\]
this implies that
\[
H(k_{\varepsilon}f)(y)=H(f)(y), \ \ \  H(1_{f(h(y))})(y)=H(k_{\varepsilon}1_{f(h(y))})(y).
\]
Thus $Hf(y)=H(1_{f(h(y))})(y)$.
(b) Suppose that $Hf(y)=H(1_{f(h(y))})(y)$ for all $f\in C(X,E)$.
Let the net $f_{\alpha}\in C(X,E)$ converges to $f$. By
Definitions \ref{d2.8} and \ref{d2.14}. Hence for big $\alpha$ we have
$\|f(h(y))-f_{\alpha}(h(y))\|\leq
\|f-f_{\alpha}\|<\delta_{f(h(y)),y}$, and so
\begin{align*}
\|Hf(y)-Hf_{\alpha}(y)\|&\leq M_y\|H(f-f_{\alpha})(y)\|\\
&=M_y\|H(1_{f(h(y))-f_{\alpha}(h(y))})(y)\|\\
&\leq M_{y}D\|H(1_{e_0})\|.\|f-f_{\alpha}\|\rightarrow 0.
\end{align*}
Therefore, $\hat{y}\circ H$ is continuous.
\end{proof}
In Theorem 5.11 of \cite{be3}, the map $H$ must be strongly
pointwise subadditive. Because, the authors used this property in
the proof. The next theorem can also be obtained with direct
computation and some minor modifications of Theorem 5.11 of
\cite{be3}, so we omit it's proof.
\begin{thm}\label{t2.19} Let $H:C(X,E)\rightarrow
C(Y,E)$ be strongly subadditive. Then the followings hold.\\
(a) If $H$ is continuous at 0, then $H$ is continuous.\\
(b) If $H$ is continuous at 0 and strongly pointwise subadditive,
then $H$ is continuous composition map.
\end{thm}
\section{Fixed point theorem}
Due to the abundant applications of fixed point theory in several disciplines such as
economics, physics, compilers and etc., many authors are interested to study the theory
of fixed point or common fixed point to different kinds of maps such that have certain
conditions. Therefore, in this section, given the importance of this theme, we introduce
the theory of fixed point to strongly suadditive maps and provide some conditions under
which the strongly subadditive map $H$ has unique fixed point and it's result will be used to
prove the existence of solution of a nonlinear integral equation. First we recall an elementary
lemma that we need it in the sequel.
\begin{lem}\label{l3.1}\cite{KHKI}
Let $\lbrace x_{n}\rbrace$ be a sequence in a metric space for which
$\Sigma_{i=1}^{\infty}d(x_{i}, x_{i+1})<\infty$. Then  $\lbrace x_{n}\rbrace$ is a
cauchy sequence.
\end{lem}
Now we give our main theorem in this section.
\begin{thm}\label{t3.2}
Let $H:C(X,E)\rightarrow C(X,E)$ be a strongly subadditive map. Suppose
the following conditions hold:
\begin{itemize}
\item[(i)] For each $f\in C(X,E)$, there exists $\epsilon > 0$ such that
$ \Vert Hf-f\Vert <\epsilon$ ;
\item[(ii)] There exists $c >0$ such that $0<cM<1$
and $\Vert Hf \Vert \leq c\Vert f\Vert$, for all $f\in C(X,E)$, where $M$ comes from strongly
subadditivity of $H$.
\end{itemize}
Then $H$ has a unique fixed point $f_{0}$ and $\lim_{n\rightarrow\infty} H^{n}f=f_{0}$.
\end{thm}
\begin{proof}
Since $H$ is strongly subadditive and condition (ii) holds, then for every $f\in C(X,E)$ we have
\begin{align*}
\Vert Hf-H(Hf)\Vert&=\Vert Hf-H^{2}f\Vert\\
&\leq M\Vert H(f-Hf)\Vert\\
&\leq cM\Vert f-Hf\Vert.
\end{align*}
Adding $\Vert f-Hf\Vert$ to both sides of the above inequality yields
\[
\Vert f-Hf\Vert+\Vert Hf-H^{2}f\Vert\leq cM\Vert f-Hf\Vert+\Vert f-Hf\Vert,
\]
hence
\[
\Vert f-Hf\Vert-cM\Vert f-Hf\Vert\leq \Vert f-Hf\Vert-\Vert Hf-H^{2}f\Vert,
\]
and so
\[
\Vert f-Hf\Vert\leq (1-cM)^{-1}\Big[\Vert f-Hf\Vert-\Vert Hf-H^{2}f\Vert \Big].
\]
Now we define the function $\varphi:C(X,E)\rightarrow R^{+}$ by
$\varphi(f)=(1-cM)^{-1}\Vert f-Hf\Vert$, for $f\in C(X,E)$. This implies that
\[
\Vert f-Hf\Vert\leq \varphi(f)-\varphi(Hf), \ \ \ f\in C(X,E).
\]
Therefore if we fix $f\in C(X,E)$ and take $m,n \in \mathbb{N}$ with $n<m$, we get that
\begin{align*}
 \Sigma_{i=n}^{m-1}\Vert H^{i}f-H^{i+1}f\Vert&\leq \varphi(H^{n}f)-\varphi (H^{m}f)\\
&<\varphi(H^{n}f).
\end{align*}
In particular by taking $n=1$ and letting $m\rightarrow \infty$ we conclude that
$\Sigma_{i=1}^{\infty}\Vert H^{i}f-H^{i+1}f\Vert \leq \varphi(Hf)<\infty$.
By Lemma \ref{l3.1}, $\{H^{n}f\}$ is a cauchy sequence. Since $C(X,E)$ is
complete, then there exists $f_{0}\in C(X,E)$ such that $\lim_{n\rightarrow\infty}H^{n}f=f_{0}$,
and since $H$ is continuous we have
\[
f_{0}=\lim_{n\rightarrow\infty}H^{n}f=\lim_{n\rightarrow\infty}H^{n+1}f=Hf_{0}.
\]
Thus $f_{0}$ is a fixed point of $H$. Now, we shall show that $f_{0}$ is the unique
fixed point. Let $g\in C(X,E)$ be another fixed point of $H$. Then
$f_{0}=\lim_{n\rightarrow\infty}H^{n}g=g$, and this completes the proof.
\end{proof}
\begin{exam}
Let $a\geq 2$ and let $z=xe^{i\theta}\in \mathbb{C}$, $-\pi<\theta<\pi$. Define
the map
\[
c_{a}:\mathbb{C}\rightarrow \mathbb{C}, \ xe^{i\theta}\longrightarrow \frac{1}{1+(2a)^{2}}xe^{ia\theta}.
\]
Since $| c_{a}(z)|=\frac{1}{1+(2a)^{2}}| z |$, it follows that
\[
\begin{array}{rl}
| c_{a}(z+z')| &=\frac{1}{1+(2a)^{2}}| z+z' |\\
&\leq \frac{1}{1+(2a)^{2}}| z |
+\frac{1}{1+(2a)^{2}}| z' |\\
&=| c_{a}(z) |+| c_{a}(z')|.
\end{array}
\]
Therefore, $c_{a}$ is subadditive. The continuity of $c_{a}$ is clear. To show that
$c_{a}$ is strongly subadditive, given any $z_{0}\in \mathbb{C}$ we must find $\epsilon_{z_{0}}>0$
and $M>0$ such that $| c_{a}(z) - c_{a}(z_{0})|\leq M | c_{a}(z-z_{0})|$
when $| z-z_{0}| <\epsilon_{z_{0}}$. We accomplish this by showing that the
function
\[
f(z,z_{0})=\frac{| c_{a}(z) - c_{a}(z_{0})|^{2}}{(1+(2a)^{2})^{2}| z-z_{0}|^{2}}
\]
is bounded in some neighborhood of $z_{0}$. Let $z=xe^{i\theta}$ and
$z_{0}=x_{0}e^{i\theta_{0}}$ observe that
\[
f(z,z_{0})=\frac{| \frac{1}{1+(2a)^{2}}xe^{ia\theta}-\frac{1}{1+(2a)^{2}}
x_{0}e^{ia\theta_{0}} |^{2}}{ (1+(2a)^{2})^{2}| xe^{i\theta}-x_{0}
e^{i\theta_{0}}|^{2}}
=
\frac{| \frac{xe^{ia\theta}}{x_{0}e^{ia\theta_{0}}}-1 |^{2}}
{| \frac{xe^{i\theta}}{x_{0}e^{i\theta_{0}}}-1 |^{2} }=f(\frac{z}{z_{0}},1).
\]
First we consider the case $z_{0}=1$, then we have
\begin{equation}\label{(3.1)}
f(z,1)=\frac{| 1-xe^{ia\theta} |^{2}}{| 1-xe^{i\theta}|^{2}}
=\frac{1+x^{2}-2xcos(a\theta)}{1+x^{2}-2xcos(\theta)}.
\end{equation}
Now we show that $f(z,1)$ is bounded from above for $z$
in a neighborhood of $1$. For all $\theta$, the following inequalities outcome from power
series expansion of $cos(\theta)$
\[
1-\frac{\theta^{2}}{2}<cos(\theta)<1-\frac{\theta^{2}}{2}+\frac{\theta^{4}}{24},\ \ \
\frac{\theta^{2}}{2}-1>-cos(\theta)>-1+\frac{\theta^{2}}{2}-\frac{\theta^{4}}{24}.
\]
This inequalities imply that
\[
f(z,1)<\frac{1+x^{2}+2x(\frac{a^{2}\theta^{2}}{2}-1)}{1+x^{2}+2x(\frac{\theta^{2}}
{2}-1-\frac{\theta^{4}}{24})}
=\frac{(1-x)^{2}+xa^{2}\theta^{2}}{(1-x)^{2}+2x(\frac{\theta^{2}}
{2}-\frac{\theta^{4}}{24})}.
\]
 For $| \theta |<\frac{1}{10}$, we have $\frac{\theta^{2}}{2}-\frac{\theta^{4}}{24}>0$
 and
 \[
f(z,1)<\frac{(1-x)^{2}+xa^{2}\theta^{2}}{(1-x)^{2}+2x(\frac{\theta^{2}}
{2}-\frac{\theta^{4}}{24})}< 1+\frac{a^{2}}{2(\frac{1}{2}-\frac{\theta^{2}}{24})}
 \]
If $| \theta |<\frac{1}{10}$, then $\frac{1}{2}-\frac{\theta^{2}}{24}>\frac{1199}{2400}$
and it follows that $f(z,1)<1+\frac{1200a^{2}}{1199}$. Thus
$\sqrt{f(z,1)}< \sqrt{1+\frac{1200a^{2}}{1199}}$.
Now suppose that $z_{0}\neq 0$ or $1$,
choose $\epsilon>0$ such that if $| \frac{z}{z_{0}}-1|<\epsilon$, then
$arg(\frac{z}{z_{0}})<\frac{1}{10}$. Since $f(z,z_{0})=f(\frac{z}{z_{0}},1)$ it follows that
$\sqrt{f(z,z_{0})}< \sqrt{1+\frac{1200a^{2}}{1199}}$, for $| z-z_{0}|<\epsilon| z_{0}|$.
After including the case $z_{0}=0$, the value $M$ previously referred to will be
\[
M=\max \Big\{1,\sqrt{1+\frac{1200a^{2}}{1199}} \Big \}=\sqrt{1+\frac{1200a^{2}}{1199}}.
\]
Also,
\begin{itemize}
\item[(i)] for each $z\in C$, there exists $\epsilon > 0$ such that
$ | c_{a}(z)-z| <\epsilon$,
\item[(ii)] for each $z=xe^{i\theta}\in C$, we have
\[
| c_{a}(z)|=\frac{1}{1+(2a)^{2}}x<\frac{1}{1+2a^{2}}x=\frac{1}{1+(2a)^{2}}| z|,
\]
\end{itemize}
if choose $c=\frac{1}{1+2a^{2}}$, then
$0<cM=\frac{1}{1+2a^{2}}\sqrt{1+\frac{1200a^{2}}{1199}}<1$.
Thus, all the conditions of the Theorem \ref{t3.2} are satisfied and $z=0$ is the
unique fixed point of $c_{a}$.
\end{exam}
Here we recall the property $P$ for a self-map. Then in Theorem \ref{t3.5} we give some
conditions under which an strongly subadditive map has $P$-property.
\begin{definition}\cite{GJ} A map $T:X\rightarrow X$ has property $P$ if
$F(T^{n})=F(T)$ for each $n\in N$, where $F(T)$ is the set of fixed points of $T$.
\end{definition}
\begin{thm}\label{t3.5}
Suppose that all the conditions of Theorem \ref{t3.2} are satisfied. Then the map $H$
has property $P$.
\end{thm}
\begin{proof}
From Theorem \ref{t3.2}, $F(H)\neq \emptyset$. Thus $F(H^{n})\neq \emptyset$
for each $n\in \mathbb{N}$. Let $n$ be a fixed positive integer greater than $1$ and
suppose that $g\in F(H^{n})$. We claim that $g\in F(H)$. Suppose on the contrary.
If $g\neq Hg$, then by conditions $(i), (ii)$ of Theorem \ref{t3.2}, and from
strongly subaddtivity of $H$ we have
\begin{align*}
\Vert H(H^{n}g)-H(H^{n-1}g)\Vert&=\Vert Hg-g\Vert\\
&\leq M\Vert H(H^{n}g-H^{n-1}g)\Vert\\
&\leq cM\Vert H^{n}g-H^{n-1}g\Vert\\
&=cM\Vert g-H^{n-1}g\Vert.
\end{align*}
Now, by letting $n\rightarrow\infty$ and Theorem \ref{t3.2}, we conclude
\[
\lim_{n\rightarrow\infty}\Vert Hg-g\Vert \leq \lim_{n\rightarrow\infty}cM\Vert g-H^{n-1}g\Vert=0,
\]
which is a contradiction. Thus $g\in F(H)$. Therefore $F(H)=F(H^{n})$ and so $H$ satisfies property $P$.
\end{proof}
\section{Application to nonlinear integral equation}
It is a truism that differential and integral equations lie at the center of mathematics, being
the inspiration of so many theoretical advances in analysis and applying to a wide range of
situations in the natural and social sciences. The term integral equation was first used by Paul
du Bois-Reymond in $1888$.\\
Recently, some researchers have been studied the problems of existence, uniqueness and other
properties of solutions of some types of nonlinear integral equations, for example, see
\cite{MoRa,AHT,NHT,HNA,HYS,HMR}.
Let us consider the following Volterra integral equation:
\begin{equation}\label{(4.1)}
x(t)=\int_{0}^{t}m(t,s)f(s,x(s))ds
\end{equation}
for all $t\in [0,T]$, where  $m(t,s)$ is real or complex valued function that are
measurable both in $t$ and $s$ on $[0,T]$, $f:[0,T]\times X\rightarrow X$ and
$x_{0}\in X$ with $X$ being a real Banach space.\\
The objective of this section is to apply Theorem \ref{t3.2}, to study the
existence and uniqueness of solution of (\ref{(4.1)}) under the conditions
in respect of the complete metric space $C([0,T],X)$ and fixed point theory.
This nonlinear integral equation will be studied under the following conditions:\\
\begin{itemize}
\item[(i)] $\sup_{0\leq t\leq T}\int_{0}^{t}|m(t,s)|ds=L_{1}<+\infty$,
\item[(ii)] for each continuous $x:[0,T]\rightarrow E$, $f(.,x(.))$ is Pettis integrable
on $[0,T]$,
\item[(iii)] for all $x,y\in C([0,T], X)$,
\[
| \int_{0}^{t}f(s,(x+y)(s))ds |\leq| \int_{0}^{t}f(s,x(s))ds|+| \int_{0}^{t}f(s,y(s))ds|, \ \ t\in[0,T],
\]
\item[(iv)] let $M>0$ and for each $x\in C([0,T],X)$ there exists $\epsilon>0$ such that
\[
| \int_{0}^{t}\Big[f(s,x(s))-f(s,y(s))\Big]ds |\leq M| \int_{0}^{t}f(s,(x-y)(s))ds|, \ \ t\in[0,T],
\]
for each $y\in C([0,T],X)$ satisfying $\parallel x-y\parallel <\epsilon$,

\item[(v)] for each $x\in C([0,T],X)$, there exists $\epsilon>0$ such that
\[
| \int_{0}^{t}m(t,s)f(s,x(s))ds-x(t) |< \epsilon, \ \ \text{for all}\ \ t\in [0,T],
\]
\item[(vi)] for each $x\in C([0,T],X)$, there exists $L_{2}>0$ such that
\[
| f(s,x(s)) |\leq L_{2}|x(s)| , \ \ \text{for all}\ \ s\in [0,T].
\]
\end{itemize}
\begin{thm}\label{t4.1} Let $X$ be a Banach space and suppose that
the assumptions (\text{i})-(vi) hold. Then the Volterra integral equation (\ref{(4.1)})
has a unique solution in $C([0,T],X)$ for which $0<L_{1}L_{2}M<1$.
\end{thm}
\begin{proof}
Define
\[
Hx(t)=\int_{0}^{t}m(t,s)f(s,x(s))ds, \ \ t\in [0,T].
\]
Then for each $x,y\in C([0,T],X)$ and by condition $(iii)$ we have
\[
\begin{array}{rl}
\Vert H(x+y)\Vert &=\Vert \int_{0}^{t}m(t,s)f(s,(x+y)(s))ds \Vert \\
&=\sup_{0\leq t\leq T}\vert \int_{0}^{t}m(t,s)f(s,(x+y)(s))ds | \\
&\leq \sup_{0\leq t\leq T}| \int_{0}^{t}m(t,s)f(s,x(s))ds | +
\sup_{0\leq t\leq T} | \int_{0}^{t}m(t,s)f(s,y(s))ds| \\
&=\Vert Hx\Vert+\Vert Hy \Vert,
\end{array}
\]
hence $H$ is subadditive. By condition $(iv)$, let $M>0$ and for each
$x\in C([0,T],X)$ there exists $\epsilon>0$ such that
\[
\begin{array}{rl}
\Vert Hx-Hy \Vert &=\Vert \int_{0}^{t}m(t,s)f(s,x(s))ds-\int_{0}^{t}m(t,s)f(s,y(s))ds \Vert \\
&=\Vert \int_{0}^{t}m(t,s)\Big[f(s,x(s))-f(s,y(s))\Big]ds \Vert \\
&=\sup_{0\leq t\leq T}\vert \int_{0}^{t}m(t,s)\Big[f(s,x(s))-f(s,y(s))\Big]ds| \\
&\leq \sup_{0\leq t\leq T}| \int_{0}^{t}m(t,s)f(s,(x-y)(s))ds| \\
&=\Vert H(x-y)\Vert,
\end{array}
\]
for each $y\in C([0,T],X)$ satisfying $\Vert x-y \Vert<\epsilon$, that is $H$
is strongly subadditive.\\
From condition $(v)$, for each $x\in C([0,T],X)$, there exists $\epsilon>0$ such that
\[
|  \int_{0}^{t}m(t,s)f(s,x(s))ds -x(t) |< \epsilon,
\]
for all $t\in [0,T]$. Then $\sup_{0\leq t\leq T}| \int_{0}^{t}m(t,s)f(s,x(s))ds- x(t)|< \epsilon,$
so $\Vert Hx-x\Vert < \epsilon$. Now, by condition $(vi)$  for each $x\in C([0,T],X)$, there exists
$L_{2}>0$ such that
\[
\begin{array}{rl}
\Vert Hx \Vert &=\Vert \int_{0}^{t}m(t,s)f(s,x(s))ds \Vert \\
&=\sup_{0\leq t\leq T}| \int_{0}^{t}m(t,s)f(s,x(s))ds| \\
&\leq \sup_{0\leq t\leq T} \int_{0}^{t}| m(t,s)| | f(s,x(s))| ds \\
&\leq \sup_{0\leq t\leq T} L_{2}\int_{0}^{t}| m(t,s)| | x(s)| ds\\
& \leq L_{1}L_{2}\Vert x\Vert.
\end{array}
\]
Thus all conditions of Theorem \ref{t3.2} are satisfied. Hence there exists a unique fixed
point $x\in C([0,T],X)$ such that $Hx=x$, which proves the existece of a uniqe solution of
(\ref{(4.1)}).
\end{proof}
Here we give an example to illustrate the usefulness of our result.
\begin{exam} In equation (\ref{(4.1)}), we define:
\[
m(t,s)=e^{t-s}sint,\ \ f(s,x)=s| x |, \ \ s,t \in[0,1], \ x\in C([0,1], \mathbb{R}),
\]
and metric $d(x,y)=\sup_{0\leq t\leq 1}| x(t)-y(t) |$ on $C([0,1], \mathbb{R})$. Then clearly
$C([0,1], \mathbb{R})$ is a complete metric space. we have
\begin{itemize}
\item[(i)] \[
\sup_{0\leq t\leq 1}\int_{0}^{t}|m(t,s)|ds=\sup_{0\leq t\leq 1}\int_{0}^{t}|e^{t-s}sint|ds
=(e-1)sin1=L_{1}<+\infty,
\]
\item[(ii)] for each continuous $x:[0,1]\rightarrow E$, $f(.,x(.))$ is Pettis integrable
on $[0,1]$,
\item[(iii)] for all $x,y\in C([0,1], \mathbb{R})$,
\[
\int_{0}^{t}f(s,(x+y)(s))ds=\int_{0}^{t}s| x(s)+y(s)| ds\leq \int_{0}^{t}s| x(s)| ds+
\int_{0}^{t}s| x(s)| ds, \ t\in [0,1],
\]
then for all $x,y\in C([0,1],\mathbb{R})$,
\[
\begin{array}{rl}
| \int_{0}^{t}f(s,(x+y)(s))ds|&=| \int_{0}^{t}s| x(s)+y(s)| ds|\leq
|\int_{0}^{t}s| x(s)| ds|+|\int_{0}^{t}s| x(s)| ds|\\
&=|\int_{0}^{t}f(s,x(s))ds|+|\int_{0}^{t}f(s,y(s))ds|, \ t\in [0,1],
\end{array}
\]
\item[(iv)] for each $x,y\in C([0,1],\mathbb{R})$ we have
\[
\begin{array}{rl}
\int_{0}^{t}\Big[f(s,x(s))-f(s,y(s))\Big]ds&=\int_{0}^{t}\Big[s | x(s)|-s | y(s)| \Big]ds\\
&\leq \int_{0}^{t}s| x(s)-y(s)| ds, \ t\in [0,1],
\end{array}
\]
then for $M=1$ and for all $x\in C([0,1],\mathbb{R})$, there exists $\epsilon>0$ such that
\[
\begin{array}{rl}
| \int_{0}^{t}\Big[f(s,x(s))-f(s,y(s))\Big]ds |&=|\int_{0}^{t}\Big[s | x(s)|-s | y(s)| \Big]ds|\\
&\leq M| \int_{0}^{t}s| x(s)-y(s)| ds|\\
&= M| \int_{0}^{t}f(s,(x-y)(s) ds|, \ t\in [0,1],
\end{array}
\]
for each $y\in C([0,1],\mathbb{R})$ satisfying $\parallel x-y\parallel <\epsilon$,

\item[(v)] for each $x\in C([0,1],\mathbb{R})$,
\[
\begin{array}{rl}
| \int_{0}^{t}m(t,s)f(s,x(s))ds-x(t) |&\leq \int_{0}^{t} |m(t,s)||f(s,x(s))|ds+
|x(t)|\\
&\leq L_{1}M_{1}+M_{1}=M_{1}(L_{1}+1),
\end{array}
\]
hence, for each $x\in C([0,1],\mathbb{R})$, there exists $\epsilon >0$ such that
\[
| \int_{0}^{t}m(t,s)f(s,x(s))ds-x(t) |\leq M_{1}(L_{1}+1)<\epsilon, \ t\in [0,1],
\]
\item[(vi)] for each $x\in C([0,1],\mathbb{R})$ and for $1=L_{2}\geq s\geq 0$,
\[
| f(s,x(s)) |=|s|x(s)||\leq L_{2}|x(s)| .
\]
\end{itemize}
Also, $0<L_{1}L_{2}M=(e-1)\sin1<1$. Thus all the conditions of Theorem \ref{t4.1} are satisfied. Hence,
the integral equation (\ref{(4.1)}) for this example, has a unique solution in $C([0,1],\mathbb{R})$.
\end{exam}


\begin{thebibliography}{99}

\bibitem{AHT}
R. P.  Agarwal, N. Hussain and M. -A. Taoudi.
Fixed point theorems in ordered Banach spaces and applications to nonlinear
integral equations. Abstract Appl. Anal. 2012, Article ID 245872, 15 pp.

\bibitem{ar3} j.Araujo E.Beckenstein and L.Narici.  Separating maps and the
non-Archimedean Hewitt theorem. Annales Math. Blaise Pascal.
2(1995) 19-27.

\bibitem{ar4} j.Araujo, E.Beckenstein and L.Narici. When is separating map
biseparating?. Archiv der Math. 67(1996)  395-407.

\bibitem{aj} j.Araujo and K. Jarosz. Automatic continuity of
biseparating maps. Studia Math. 155 (2003), no. 3, 231-239.

\bibitem{be3} E.Beckenstein and L.Narici. Subadditive separating maps Acta math
Hungar. 88(1-2) 2000  147-167.

\bibitem{be2} J. J. Font and S. Hernandez. On separating maps between locally
compact spaces. Arch. Math. 63(1994) 158-165.

\bibitem{NHT}
N. Hussain and M.-A. Taoudi.
Krasnosel'skii-type fixed point theorems with applications to Volterra integral equations.
Fixed Point Theory Appl., 2013, 2013:196.

\bibitem{ja} K. Jarosz. Automatic continuity of separating linear
isomorphism. Canad. Math. Bull. 33 (1990), no. 2, 139-144.

\bibitem{GJ}
G. S. Jeong and B. E. Rhodes.
More maps for which $F(T^{n})=F(T)$. Demonstratio Mathematica Vol. XL, no.  3 (2007), 671-680.

\bibitem{MoRa}
B. Moeini and A. Razani.
$\mathcal{JH}$-operator pairs of type $(R)$ with application to nonlinear integral equations.
Vietnam journal of mathematics (accepted).

\bibitem{nai1} M .Naimark. Normed Rings.  Nordhoff (Groningen 1959).

\bibitem{HNA}
H. K. Pathak, S. N. Mishra and A. K. Kalinede and S.S. Chang.
Common fixed point theorems with applications to nonlinear integral equations.
Demonstratio Math. XXXII (3) (1999), 547-564.

\bibitem{HYS}
H. K. Pathak, Y. J. Cho and S. M. Kang.
Common fixed points of biased maps of type $(A)$ and applications. Internat. J. Math.
Math. Sci. 21 (1998), no. 4, 681-694.

\bibitem{HMR}
H. K. Pathak, M. S. Khan and Rakesh Tiwari.
Common fixed point theorem and application to nonlinear integral equations.
Computers and Mathematics with Applications 53 (2007), 961-971.

\bibitem{KHKI}M. A. Khamsi and W. A. Kirk.
An introduction to metric spaces and fixed point theory. Pure and Applied Mathematics:
A wiley-Interscience series of Texts, Monographs, and Tracts (2001).

\bibitem{se} F. Sady Y. Estaremi. Subadditive separating maps between regular Banach
function algebras. Bull. Korean Math. Soc. 44 (2007)  753-761.

\bibitem{vv} V. Valvo, On separating subadditive maps. Turk. J. Math 39 (2015)  168-173.

\end{thebibliography}
\end{document}